\documentclass[preprint,3p,12pt,pdf]{elsarticle}

\makeatletter
\def\ps@pprintTitle{%
  \let\@oddhead\@empty
  \let\@evenhead\@empty
  \def\@oddfoot{\reset@font\hfil\thepage\hfil}
  \let\@evenfoot\@oddfoot
}
\makeatother

\usepackage{lineno,hyperref}
\modulolinenumbers[5]

\usepackage{hyperref}

\usepackage{epstopdf}

\usepackage{tikz}
\usetikzlibrary{arrows}

\usepackage{pst-node}
\usepackage{tikz-cd}

\usepackage[shellescape]{gmp}

\usepackage{mathrsfs}
\usepackage{amssymb}
\usepackage{amsfonts}
\usepackage{latexsym}
\usepackage{mathtools}
\usepackage{xcolor}
\usepackage{wrapfig}
\usepackage{floatflt}

\usepackage{mathtools}
\usepackage{extarrows}

\usepackage{graphicx}

\usepackage{subcaption}

\usepackage{epigraph}

\def\lb{\label}

\newcommand{\er}[1]{\textrm{(\ref{#1})}}

\newtheorem{theorem}{\bf Theorem}[section]

\def\a{\alpha}

\def\d{\delta}

\def\z{\zeta}          
           
\def\p{\psi}

\def\m{\mu}

       \def\vp{\varphi}    

       \def\C{{\mathbb C}}    
    \def\N{{\mathbb N}}   

\def\lt{\biggl}                  \def\rt{\biggr}
\def\ol{\overline}               \def\wt{\widetilde}

\let\ge\geqslant                 

\def\iy{\infty}

\def\el2{\ell^{\,2}}             \def\1{1\!\!1}

\def\Im{\mathop{\mathrm{Im}}\nolimits}

\def\Re{\mathop{\mathrm{Re}}\nolimits}

\newtheorem{corollary}[theorem]{\bf Corollary}

\let\ge\geqslant

\newcommand{\ca}{\begin{cases}}
	\newcommand{\ac}{\end{cases}}
\newcommand{\ma}{\begin{pmatrix}}
	\newcommand{\am}{\end{pmatrix}}
\renewcommand{\[}{\begin{equation}}
	\renewcommand{\]}{\end{equation}}
\def\eq{\begin{equation}}
	\def\qe{\end{equation}}

\begin{document}
	
	\begin{frontmatter}

		\title{Closed-form solutions for Bernoulli and compound Poisson branching processes in random environments}

		\date{\today}

		\author
		{Anton A. Kutsenko}
	
		\address{University of Hamburg, MIN Faculty, Department of Mathematics, 20146 Hamburg, Germany; email: akucenko@gmail.com}
		\address{Mathematical Institute for Machine Learning and Data Science, Katholische Universit\"at Eichst\"att--Ingolstadt, Germany\footnote{The author was affiliated with the Katholische Universit\"at Eichst\"att--Ingolstadt at the time of the research initiation and is currently affiliated with the University of Hamburg.}}

\epigraph{In Memory of Pavel Petrovitch Kargaev (1953-2024)}	
	
	\begin{abstract}
		For branching processes, the generating functions for limit distributions of so-called ratios of probabilities of rare events satisfy the Schr\"oder-type integral-functional equations. Excepting limited special cases, the corresponding equations can not be solved analytically. I found a large class of Poisson-type offspring distributions, for which the Schr\"oder-type functional equations can be solved analytically. Moreover, for the asymptotics of limit distributions, the power and constant factor can be written explicitly. As a bonus, Bernoulli branching processes in random environments are treated. The beauty of this example is that the explicit formula for the generating function is unknown. Still, the closed-form expression for the power and constant factor in the asymptotic can be written with the help of some outstanding tricks. Also, the asymptotic expansion contains oscillatory terms absent in the Poisson case. It is shown that at the power $3$ in the Bernoulli binomial kernels short-phase discrete oscillations turn into long-phase continuous ones.
	\end{abstract}


	
\end{frontmatter}


{\section{Introduction}\lb{sec0}}

One of the motivations for this work is a limited number of closed-form solutions of Schr\"oder-type equations related to branching processes. I will add a few precious stones to this piggy bank.

Let us start with general results related to branching processes in random environments and limit distributions of rare events. Let $X_t$ be a Galton-Watson branching process in random environments with offspring probability-generating function 
\[\lb{001}
 P_r(z)=p_{1r}z+p_{2r}z^2+p_{3r}z^3+....
\]
Here $p_{jr}$ are the probabilities of producing $1$, $2$, ... offspring by individuals in a generation - they are all non-negative and their sum is $\sum_{j=1}^{+\iy} p_{jr}=1$. At each time $t$, the parameter $r$ is chosen randomly from some probabilistic space $(R,\m)$. The function $P_r(z)$ is the same for all individuals at fixed time $t$. At the next time $t+1$, the function $P_{\wt r}(z)$, again the same for all individuals at this time $t+1$, can differ from $P_r(z)$ because $\wt r$ is randomly chosen. 

Since the minimal family size is at least one ($p_{0r}=0$, and $p_{1r}\ne0$ is assumed) in each generation, one can define the so-called {\it limit of ratios of probabilities of rare events}
\[\lb{002}
 \vp_n:=\lim_{t\to+\iy}\frac{\mathbb{P}(X_t=n)}{\mathbb{P}(X_t=1)},
\]
and the corresponding generating function
\[\lb{003}
 \Phi(z)=\vp_1z+\vp_2z^2+\vp_3z^3+....
\]
This function satisfies the Schr\"oder-type integral-functional equation
\[\lb{004}
 \int_R\Phi(P_r(z))\m(dr)=\Phi(z)\int_Rp_{1r}\m(dr),\ \ \ \Phi(0)=0,\ \ \ \Phi'(0)=1.
\] 
A short justification of \er{004} is given in the Appendix. The details, including the existence of solution of \er{004}, some interesting properties, and examples are provided in \cite{K1}. 
Note that, since $\int\m(dr)$ appears on both sides of \er{004}, it is not necessary to assume $\m(R)=1$, it is sufficient to have $\m(R)>0$ plus the convergence of integrals. The difference between \er{004} and the classical Schr\"oder functional equation is the presence of these integrals. Integral operators in LHS of \er{004} are generally more complex than known types of operators, e.g., Hilbert-Schmidt or Calder\'on-Zygmund, see for example \cite{CGPSZ}, and they also may lead to various types of oscillations in the solutions of corresponding integral equations, see Examples below. The oscillations, or the lack thereof, are the most interesting subjects related to these equations.

Let us discuss briefly the physical meaning and motivation for the study of the ratios $\vp_n$ in \er{002}. Since $p_{0r}=0$, the branching process grows. However, 
there is a non-zero probability at each time step that only one particle will survive. 
It is natural to compare this probability with the probabilities that only $2$, $3$, ... particles will survive. The formula \er{002} answers this question. The mentioned events are rare, but their probability increases with increasing $n$. The asymptotic of $\vp_n$ for large $n$ is of special interest. For the classical Galton-Watson process, this asymptotic is closely related to the left tail asymptotic of the density $p(x)$ of the {\it martingal limit} $\lim_{t\to+\iy}E^{-t}X_t$, where $E$ is the expectation. This density shows how the population grows overall, excluding rare events. Researchers pay the greatest attention to tail asymptotics of $p(x)$, when $x\to+0$ and $x\to+\iy$. The analysis of the behavior of $p(x)$ at the right tail, $x\to+\iy$, was started in \cite{H}, where oscillations were noted, and then continued in \cite{BB,FW}, and \cite{FW1}. Although the asymptotic series obtained are incomplete, the analysis appears to have been finished more or less. The left tail differs from the right tail drastically. Significant progress in describing the structure of the first asymptotic term of $p(x)$ on the left tail, when $x\to+0$, was made in \cite{BB}. In the recent paper \cite{K24}, the complete asymptotic series for $p(x)$ was obtained. Moreover, the series converges for all $x>0$, not only for small $x\to+0$. One of the key ingredients for obtaining the series was that the first asymptotic term of $\vp_n$ for $n\to+\iy$ and of $p(x)$ for $x\to+0$ are, in some sense, the same. I believe that the analysis methods of $\vp_n$ will also help us study $p(x)$ for the more complex case of random environments. At the moment, these questions are largely open.

{\bf Remark on explicit solutions.} As in the classical case of the constant environment, except, e.g., the geometric distributions $P_r(z)=rz/(1-(1-r)z)$, equation \er{004} can not be solved analytically in general. At least, I do not know the corresponding methods. The geometric distribution is a partial case of the following class of probability-generating functions
\[\lb{expl}
 P_r(z):=P^{-1}(rP(z)).
\]
For \er{expl}, the classical Schr\"oder equation, as well as the case of random environments by the parameter $r$, can be solved explicitly, namely $\Phi(z)=P(z)$. The case $P(z)=z/(1-z)$ corresponds to the geometric distribution. Another interesting example $P(z)=-\ln(1-z)$  leads to $P_r(z)=1-(1-z)^{r}$. In both cases, $r\in(0,1]$ - otherwise, $P_r(z)$ is not a probability generating function. Regarding \er{expl} and possible generalizations, see, e.g., \cite{SL2016} and \cite{GH2017}. 

Apart from the case \er{expl}, I found another large one-parametric family that can be treated explicitly. Let
\[\lb{e005}
 q(z)=q_1z+q_2z^2+q_3z^3+...,\ \ \ all\ \ q_j\ge0
\]
be some analytic function in the unit disc satisfying $0<q(1)<+\iy$. Then, for any $r\in(0,+\iy)$,
\[\lb{006}
 P_r(z)=z\frac{e^{r q(z)}}{e^{rq(1)}} 
\] 
is a probability-generating function of the form \er{001}. The significant difference between \er{expl} and \er{006} is that no explicit solutions of classical Schr\"oder functional equations related to \er{006} are known. Such probability-generating functions are directly related to {\it Compound Poisson Distributions}. Indeed, consider
\[\lb{CPD}
\xi=1+\sum_{n=1}^N \eta_n,\ \ \ N\sim {\rm Poisson}(rq(1)),
\] 
where $\eta_n$ are random variables, with the same probability generating function $q(z)/q(1)$, that are mutually independent and also independent of $N$. Then, due to the {\it law of total expectation}, see, e.g. \cite{H}, the probability generating function of $\xi$ is given by \er{006}. If $q(z)=z$ then $P_r(z)$ corresponds to a simple non-degenerating Poisson distribution. We assume that $r$ is taken randomly and uniformly from $(0,+\iy)$ - i.e., the corresponding probability measure is $({\rm Lebesgue\ measure})/R$ on $[0,R)$ and then we take the limit $R\to+\iy$. This means that the corresponding Galton-Watson process mentioned at the beginning of the Introduction is given by
\[\lb{GWP}
X_{t+1}=\sum_{j=1}^{X_t}\xi_{j,t},\ \ \ X_0=1,\ \ \ t\in\N\cup\{0\},
\]
where $\xi_{j,t}$ are all independent and have the form \er{CPD} with $r=r_t$ taken randomly and uniformly from $(0,+\iy)$, and independently at each time step $t$.  Stochastic processes similar to compound Poisson branching ones are used in practice, see, e.g., \cite{BJR,P}, and \cite{CH}.  

Substituting \er{006} into \er{004} and taking into account the fact that $p_{1r}=e^{-rq(1)}$, we obtain
\[\lb{007}
 \int_{0}^{+\iy}\Phi(z{e^{r (q(z)-q(1))}})dr=\frac{\Phi(z)}{q(1)},\ \ \ \Phi(0)=0,\ \ \ \Phi'(0)=1.
\]
The result will be the same for $q(z)$ and $q(z)/q(1)$. Hence, without loss of generality, one may consider $q(1)=1$, but we proceed without this assumption. Equation \er{007} can be solved explicitly.
\begin{theorem}\lb{T1}
	The solution of \er{007} is
\[\lb{008}
 \Phi(z)=z\exp\lt(\int_0^z\frac{ q'(\zeta)\zeta+q(\zeta)}{(q(1)-q(\z))\zeta}d\zeta\rt).
\]
\end{theorem}

For polynomials $q(z)$, further simplification of \er{008} is possible. We formulate one in the simplest case.

\begin{corollary}\lb{C1}
Suppose that $q(z)$ is a polynomial of degree $N\ge1$. Let $b_j$, $j=1$,...,$N$ be all the zeros of equation $q(z)=q(1)$. If they are all distinct then
\[\lb{009}
 \Phi(z)=z\prod_{j=1}^N\lt(1-\frac{z}{b_j}\rt)^{-1-\frac{q(1)}{q'(b_j)b_j}}.
\]
\end{corollary}

If multiple roots exist, then \er{009} requires small modification. Anyway, having such type of formulas at your disposal, one may find explicit expressions of $\vp_n$ through extended binomial coefficients
\[\lb{bin}
 \binom{\a}{0}=1,\ \ \ \binom{\a}{n}=\frac{\a\cdot...\cdot(\a-n+1)}{1\cdot...\cdot n},\ \ n\ge1,
\]
and compute asymptotics of $\vp_n$ for large $n$ explicitly.

\begin{corollary}\lb{C2}
Under the assumptions of Corollary \ref{C1}, we have
\[\lb{010}
 \vp_n=\sum_{\sum_{j=1}^Nm_j=n-1,\ m_j\ge0}\prod_{j=1}^N\binom{-1-\frac{q(1)}{q'(b_j)b_j}}{m_j}\lt(\frac{-1}{b_j}\rt)^{m_j}.
\]
Moreover, if there is one only among $b_j$ equal to $1$, say $b_1=1$, and all others satisfy $|b_j|>1$ then, for large $n$, the following asymptotic expansion is fulfilled
\[\lb{011}
 \vp_n\simeq C\lt(n^{\frac{q(1)}{q'(1)}}+\frac{q(1)(q'(1)-q(1))(q''(1)+q'(1))}{2q'(1)^3}n^{\frac{q(1)}{q'(1)}-1}+...\rt),
\]
where
\[\lb{012} 
 C=\frac{1}{\Gamma(1+\frac{q(1)}{q'(1)})}\prod_{j=2}^N\lt(1-\frac{1}{b_j}\rt)^{-1-\frac{q(1)}{q'(b_j)b_j}}.
\]  
\end{corollary}

The method presented in the Proof Section allows one to obtain the complete asymptotic series. The second term in \er{011} is already attenuating. At the moment, I have no good estimates for the remainder. Nevertheless, as the examples below show, even two terms are almost perfect. All the asymptotic terms do not contain any oscillatory factors. Typically, oscillations appear in the first and other terms for the classical case of the Galton-Watson process. For the Galton-Watson processes in random environments, they usually appear not in the first but in the next asymptotic terms, see \cite{K1}. Oscillations in the limit distributions of branching processes and their absence may be important in practical applications, see, e.g., \cite{DIL,CG}, and \cite{DMZ}. Roughly speaking, one naturally expects something like power-law growth of ratios of probabilities of rare events $\vp_n$ when $n$ grows, but, often, periodic factors should correct this growth. The exact form of the correction factors can be then used to predict or model necessary effects: increase or decrease the amplitudes of oscillations. The absence of oscillations in the leading asymptotic term of ratios of probabilities of rare events may indicate that the system lives in random environments, and standard (constant environment) Galton-Watson processes may be inappropriate for modeling its evolution.

We have obtained explicit solutions of the generalized Schr\"oder-type integral-functional equation for a large class of Poisson-type kernels. For the classical case of the Galton-Watson process, the standard Schr\"oder functional equation usually cannot be solved explicitly. One of the reasons is that the domain of definition for the solution is related to the filled Julia set of the probability-generating function. Very often, this set has a complex fractal structure. For example, for the Poisson-type exponential probability generating functions considered above, this is a so-called ``Cantor bouquet'', see Example 3 below. Explicit functions with such a complex domain of definition are practically unknown. The same can be said about the asymptotics of the power series coefficients of these solutions. It is only known that the asymptotic consists of periodic terms multiplied by some powers, see, e.g., \cite{B1,BB}, and \cite{K}. Generally speaking, there are no explicit representations for the corresponding Fourier coefficients of the periodic terms through known constants. Galton-Watson processes in random environments are similar to the classical case in that there are usually no explicit expressions for the solutions of the generalized Schr\"oder-type integral-functional equations. However, there are some differences. For example, the periodic terms usually disappear in the leading asymptotic term of the limit distribution, see also Example 4 below. They may appear in the next asymptotic terms, but, in contrast to the classical case, the fast algorithms for their computations are usually unknown, see \cite{K1} and \cite{K}. Sometimes, with a big effort leading asymptotic terms can be computed explicitly. For example, if $P_r(z)=rz+(1-r)z^2$ with uniform $r\in(0,1)$ then
\[\lb{onetwof}
 \lim_{t\to+\iy}\frac{\mathbb{P}(X_t=n)}{\mathbb{P}(X_t=1)}=\frac{n}{2-\ln 4}+...,\ \ \ n\to+\iy.
\] 
The spectacular story about the constant $(2-\ln4)^{-1}$ is presented in \cite{K1}. Apart from the already mentioned case \er{expl}, such explicit results are practically unknown. We formulate an extension of \er{onetwof} to binomial kernels of arbitrary power.
\begin{theorem}
\lb{T2}
 If $P_r(z)=z(r+(1-r)z)^m$ with uniform $r\in(0,1)$ then
\[\lb{T2f1}
 \lim_{t\to+\iy}\frac{\mathbb{P}(X_t=n)}{\mathbb{P}(X_t=1)}=\frac{m^2}{(m+1)(m-\ln(m+1))}(n+\p_n+...),\ \ \ n\to+\iy,
\]
where $m+1$-periodic piece-wise linear sequence $\p_n$, $n\ge1$, is defined by
\[\lb{T2f2}
 \p_n=\ca 
 \p_{n-1}+\frac1m,& m+1\nmid n-1,\\
 \p_{n-1}-1,& otherwise
 \ac,\ \ \ \p_{m+1}=\frac{m+1}{2m}\cdot\frac{(m+2)\ln(m+1)-2m}{(m+1)\ln(m+1)-m}.
\]
\end{theorem}
The proof of Theorem \ref{T2} is given in Example 5 below. Asymptotic \er{T2f1} should be interpreted as an approximation because we do not specify the growth rate for the reminder. 
However, examples show that this approximation is almost perfect. Interestingly, the third asymptotic term is a sequence of period $(m+1)^2$ multiplied by $n^{-1}$ for the first three values $m=1,2,3$. Starting from $m=4$, this is a long-phase oscillation multiplied by $n^{-\a}$ with $\a\in(0,1)$ is a first root of \er{304b} with a positive real part.

Summarizing the discussion, I have found an extensive class of Poisson-type distributions, where everything can be calculated explicitly. It is important to note that, in this class, there are no oscillations in the asymptotic terms. However, in almost all other examples, periodic terms are presented in the asymptotic. Even in \er{onetwof}, there are infinitely many non-leading oscillatory asymptotic terms.

The rest of the story contains proof of the main results and examples. The first three Examples (0, 1, and 2) below directly relate to the main results. Examples 3 and 4 are not a part of the main results but are extremely useful for comparison of the cases of random environments in complete and partial parametric intervals and constant environments. There are two main differences here: the presence of oscillations, and the absence of explicit formulas for multipliers in the asymptotic. The bonus example 5 is wonderful - it contains short and long phase oscillations, and the multipliers here can be computed with the help of some intriguing tricks.

The Appendix contains little hints on how the Schr\"oder equation is derived.

\section{Proof of Theorem \ref{T1} and Corollaries \ref{C1} and \ref{C2}}
The main result \er{008} follows from the following reasoning based on the change of variables
\[\lb{100a}
 \int_{0}^{+\iy}\Phi(z{e^{r (q(z)-q(1))}})dr=\lt[w=z{e^{r (q(z)-q(1))}},\ dr=\frac{dw}{w(q(z)-q(1))}\rt]=
\]
$$
\frac1{q(1)-q(z)}\int_0^z\frac{\Phi(w)}{w}dw,\ \ \ |z|<1,
$$
where the fact that $\Re(q(z)-q(1))<0$ inside the unit disc is used for obtaining $w(r=+\iy)=0$ as one of the integration limits. Introducing $\wt\Phi(z)=\Phi(z)/z$, and using \er{007}, we obtain the simple ODE
$$
 \wt\Phi(z)=\frac{(z(q(1)-q(z))\wt\Phi)'(z)}{q(1)},\ \ \ \wt\Phi(0)=1,
$$
solution of which is
$$
 \wt\Phi(z)=\exp\lt(\int_0^z\frac{ q'(\zeta)\zeta+q(\zeta)}{(q(1)-q(\z))\zeta}d\zeta\rt),
$$
which, at the end, gives \er{008}.

If all the zeros $\{b_j\}_{j=1}^N$ of the polynomial $q(1)-q(\zeta)$ are distinct then the fraction decomposition 
$$
 \frac{q'(\zeta)\zeta+q(\zeta)}{(q(1)-q(\z))\zeta}=\sum_{j=1}^N\frac{1+\frac{q(b_j)}{q'(b_j)b_j}}{b_j-\zeta}
$$
along with the identities $q(b_j)=1$ and \er{008} give \er{009}. Identity
\[\lb{100}
 (1+w)^{\a}=\sum_{j=0}^{+\iy}\binom{\a}{j}w^j
\]
along with \er{009} give \er{010}. Under the assumptions of Corollary \ref{C1}, we see that the main contribution to the asymptotic comes from the smallest singularity of $\Phi(z)$ - from the point $z=1$. It is appropriate here to mention the work \cite{FO}, which contains many remarkable techniques of singularity analysis and transfer theorems. Some of the ideas from this work will continue to appear in one form or another. Due to \er{009} and definition \er{012}, the main term is
$$
 \Phi(z)=C\Gamma\lt(1+\frac{q(1)}{q'(1)}\rt)(1-z)^{-1-\frac{q(1)}{q'(1)}}+...,
$$ 
which, by \er{100}, leads to
\[\lb{101}
 \vp_n\simeq C\Gamma(1+\frac{q(1)}{q'(1)})(-1)^n\binom{-1-\frac{q(1)}{q'(1)}}{n}+....
\]
Using \er{101} and the known asymptotic expansion of generalized binomial coefficients, see \cite{TE} and \cite{K},
\[\lb{101a}
 \binom{\a}{n}\simeq\frac{(-1)^n}{\Gamma(-\a)n^{\a+1}}+\frac{(-1)^n\a(\a+1)}{2\Gamma(-\a)n^{\a+2}}+...,
\]
we obtain the main term in \er{011}. In principle, using \er{009}, we can continue with expansions of $\Phi(z)$ at $z=1$
\[\lb{102}
 \Phi(z)= C\Gamma(1+\frac{q(1)}{q'(1)})((1-z)^{-1-\frac{q(1)}{q'(1)}}+A_1(1-z)^{-\frac{q(1)}{q'(1)}}+A_2(1-z)^{1-\frac{q(1)}{q'(1)}}+...)
\]
to obtain $A_j$, and then all the asymptotic terms of $\vp_n$. However, there is a more fast, than \er{009}, algorithm for obtaining $A_j$ based on the integral equation \er{007} and change-of-variables presented at the beginning of this section. Using \er{100a} and \er{004}, we have 
$$
 \int_0^z\frac{\Phi(w)}{w}dw=(q(1)-q(z))\frac{\Phi(z)}{q(1)}
$$
or, by \er{102} and the geometric progression formula,
$$
 \int_0^zC\Gamma(1+\frac{q(1)}{q'(1)})((1-w)^{-1-\frac{q(1)}{q'(1)}}+A_1(1-w)^{-\frac{q(1)}{q'(1)}}+...)(1+(1-w)+(1-w)^2+...)dw\simeq
$$
$$
 (q'(1)(1-z)-\frac{q''(1)}2(1-z)^2+...)\frac{C\Gamma(1+\frac{q(1)}{q'(1)})}{q(1)}((1-z)^{-1-\frac{q(1)}{q'(1)}}+A_1(1-z)^{-\frac{q(1)}{q'(1)}}+...).
$$
Integrating and equating the coefficients related to $(1-z)^{1-\frac{q(1)}{q'(1)}}$, we obtain
$$
 \frac{1}{\frac{q(1)}{q'(1)}-1}+\frac{A_1}{\frac{q(1)}{q'(1)}-1}=-\frac{q''(1)}{2q(1)}+\frac{A_1q'(1)}{q(1)},
$$
which gives
\[\lb{103}
 A_1=\frac{(q'(1)-q(1))q''(1)-2q(1)q'(1)}{2(q'(1))^2}.
\]
Using \er{102} and \er{101a} with \er{103}, we obtain the explicit expression for the second term in \er{011}:
$$
 \frac{(q(1)+q'(1))q(1)}{2(q'(1))^2}+\frac{\Gamma(1+\frac{q(1)}{q'(1)})}{\Gamma(\frac{q(1)}{q'(1)})}A_1=
 \frac{q(1)(q'(1)-q(1))(q''(1)+q'(1))}{2(q'(1))^3}.
$$

\section{Examples \lb{sec2}}

For most examples, Embarcadero Delphi Rad Studio Community Edition and the library NesLib.Multiprecision is used. This software provides a convenient environment for programming and well-functioning basic functions for high-precision computations. All the algorithms related to the article's subject, including efficient parallelization, are developed by the author (AK).  

\subsection{Example 0.} The case $q(z)=z$ means averaging over all the simple Poisson branching processes. We have $\Phi(z)=z(1-z)^{-2}$, see \er{009}, and, hence,
$$
 \lim_{t\to+\iy}\frac{\mathbb{P}(X_t=n)}{\mathbb{P}(X_t=1)}=n.
$$ 
In addition to geometric distributions, this is the first example of a complete closed-form solution, including simple asymptotic behavior.

\subsection{Example 1.} 
\begin{figure}
	\centering
	\begin{subfigure}[b]{0.45\textwidth}
		\includegraphics[width=\textwidth]{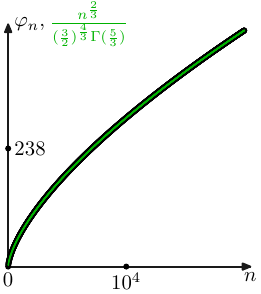}
		\caption{First asymptotic term}
		\label{fig1a}
	\end{subfigure}
	~ 
	\begin{subfigure}[b]{0.45\textwidth}
		\includegraphics[width=\textwidth]{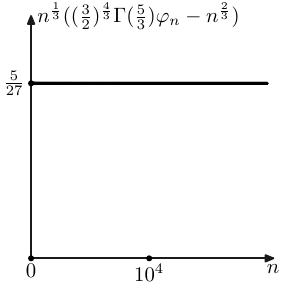}
		\caption{Second asymptotic term}
		\label{fig1b}
	\end{subfigure}
	\caption{Plots show how good asymptotic \er{200} is.}\label{fig1}
\end{figure}

The case $q(z)=z+z^2$. We have
$$
 b_1=1,\ \ b_2=-2,\ \ q(1)=q(b_1)=2,\ \ q'(1)=q'(b_1)=3,\ \ q'(b_2)=-3,\ \ q''(1)=2.
$$
The explicit expression for $\vp_n$, see \er{010}, is
$$
 \vp_n=\sum_{m_1=0}^{n-1}\binom{-\frac{5}{3}}{m_1}(-1)^{m_1}\binom{-\frac{4}{3}}{n-1-m_1}\lt(\frac{1}{2}\rt)^{n-1-m_1}.
$$
The asymptotic for large $n$, see \er{011}, is
\[\lb{200}
 \vp_n\simeq\frac{1}{(\frac32)^{\frac43}\Gamma(\frac53)}\lt(n^{\frac23}+\frac{5}{27}n^{-\frac13}+...\rt).
\]
As seen in Fig. \ref{fig1}, this asymptotic expansion is perfect enough.

\subsection{Example 2.} 
\begin{figure}
	\centering
	\begin{subfigure}[b]{0.45\textwidth}
		\includegraphics[width=\textwidth]{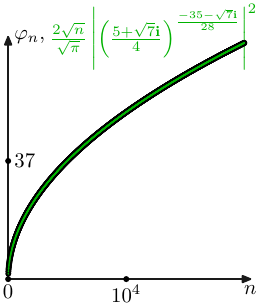}
		\caption{First asymptotic term}
		\label{fig2a}
	\end{subfigure}
	~ 
	\begin{subfigure}[b]{0.45\textwidth}
		\includegraphics[width=\textwidth]{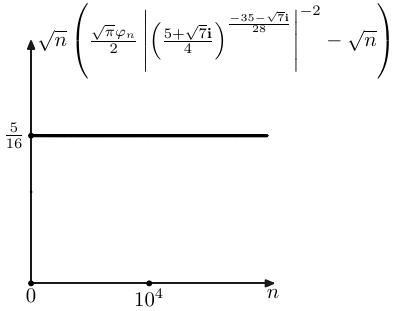}
		\caption{Second asymptotic term}
		\label{fig2b}
	\end{subfigure}
	\caption{Plots show how good asymptotic \er{201} is.}\label{fig2}
\end{figure}

The case $q(z)=z+z^3$. We have
$$
b_1=1,\ \ b_2=\frac{-1-\sqrt{7}\mathbf{i}}2,\ \ b_3=\frac{-1+\sqrt{7}\mathbf{i}}2,\ \ q(1)=q(b_1)=2,\ \ q'(1)=q'(b_1)=4,\ \ q''(1)=6,
$$
$$
 q'(b_2)=\frac{-7+3\sqrt{7}\mathbf{i}}2,\ \ q'(b_3)=\frac{-7-3\sqrt{7}\mathbf{i}}2.
$$
The explicit expression for $\vp_n$, see \er{010}, is
$$
\vp_n=\sum_{m_2=0}^{n-1}\sum_{m_1=m_2}^{n-1}\binom{-\frac{4}{3}}{m_2}(-1)^{m_2}\binom{\frac{-35+\sqrt{7}\mathbf{i}}{28}}{m_1-m_2}\lt(\frac{-1+\sqrt{7}\mathbf{i}}{4}\rt)^{m_1-m_2}\binom{\frac{-35-\sqrt{7}\mathbf{i}}{28}}{n-1-m_1}\lt(\frac{-1-\sqrt{7}\mathbf{i}}{4}\rt)^{n-1-m_1}.
$$
The asymptotic for large $n$, see \er{011}, is
\[\lb{201}
\vp_n\simeq\frac{2}{\sqrt{\pi}}\lt|\lt(\frac{5+\sqrt{7}\mathbf{i}}{4}\rt)^{\frac{-35-\sqrt{7}\mathbf{i}}{28}}\rt|^2\lt(\sqrt{n}+\frac{5}{16\sqrt{n}}+...\rt).
\]
Again, this asymptotic expansion is perfect as seen in Fig. \ref{fig2}.

\subsection{Example 3.} 
\begin{figure}
	\centering
	\begin{subfigure}[b]{0.45\textwidth}
		\includegraphics[width=\textwidth]{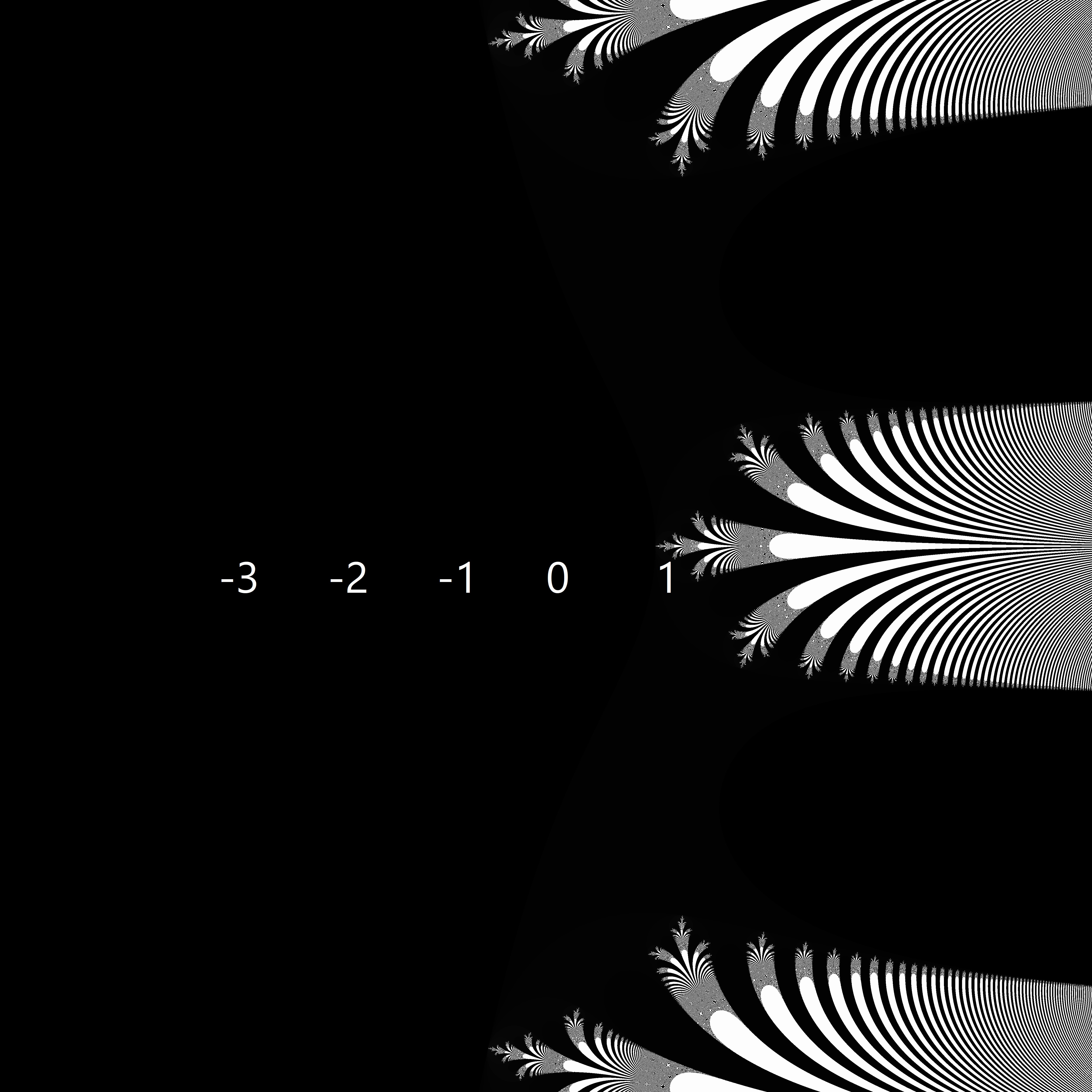}
		\caption{Filled Julia set (black area) for $ze^{z-1}$}
		\label{fig3a}
	\end{subfigure}
	~ 
	\begin{subfigure}[b]{0.45\textwidth}
		\includegraphics[width=\textwidth]{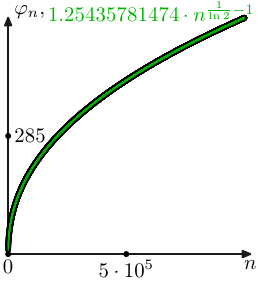}
		\caption{Exact values and first asymptotic term}
		\label{fig3b}
	\end{subfigure}
	\caption{(a) The domain of definition (black area) of $\Phi(z)$, see \er{202}, in complex plane $\C$; (b) Taylor coefficients of $\Phi(z)$ compared with its first asymptotic term.}\label{fig3}
\end{figure}

\begin{figure}
		\includegraphics[width=\textwidth]{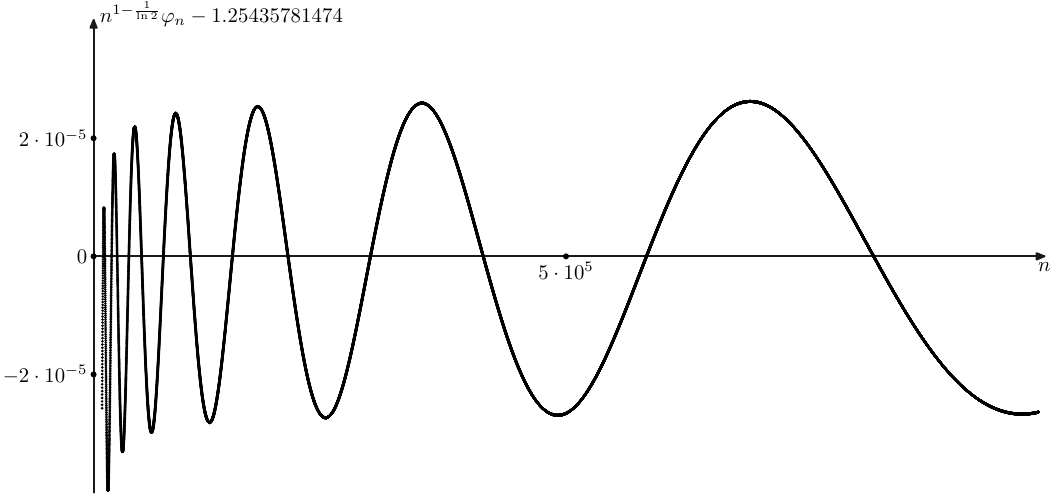}
		\caption{Oscillatory factors in the first asymptotic term for $\vp_n$, see \er{203}.}
		\label{fig4}
\end{figure}

Let us compare the cases of random environments with a Galton-Watson process in a constant environment. We take $q(z)=z$ considered in Example 0 and fix $r=1$. In other words, we take the delta measure $\d(r=1)$ in \er{004}. Then the generating function for the ratios of rare events $\Phi(z)$ satisfies the Schr\"oder-type functional equation
\[\lb{202}
 \Phi(ze^{z-1})=e^{-1}\Phi(z),\ \ \Phi(0)=0,\ \ \Phi'(0)=1,
\]
see \er{007}, where $r=1$ is fixed. Now, we will follow the ideas of \cite{K} and mostly skip some details. No explicit solutions for \er{202} are known. It is not surprising, since the domain of definition for $\Phi(z)$ coincides with the filled Julia set for the exponential generating function and has a complex fractal form, see its approximate structure in Fig. \ref{fig3a}.
Using \er{202}, we can write recurrence expression for determining $\vp_n$ - Taylor coefficients of $\Phi(z)$ at $z=0$:
$$
 \vp_1=1,\ \ \ \frac{\vp_n}{e}=\frac{\vp_n}{e^n}+\frac{(n-1)^1\vp_{n-1}}{1!e^{n-1}}+\frac{(n-2)^2\vp_{n-2}}{2!e^{n-2}}+...+\frac{\vp_{1}}{(n-1)!e},\ n>1.
$$ 
Non-recurrence formulas like \er{010} for $\vp_n$ in this example are unknown. Results of \cite{K} allow us to write an approximate expression
\[\lb{203}
\vp_n\approx1.25435781474... n^{{\frac1{\ln2}}-1},
\]
where the constant is an average value of the so-called Karlin-MacGregor function related to $ze^{z-1}$, the approximation is quite good, see Fig. \ref{fig3b}. However, in contrast to the examples of random environments, no explicit expressions through known constants exist. Moreover, apart from the constant term multiplied by the power of $n$, the first asymptotic term also contains small oscillatory elements, see Fig. \ref{fig4}. The presence of oscillations is the significant difference between constant and random environments. Note that, if we take averaging not over the entire interval $(0,+\iy)$ in \er{007}, but only over part of it, then the oscillations in the first asymptotic term will disappear, but will remain in the following terms of the asymptotic, as it is discussed in \cite{K1}. We show this explicitly in the next example.

\subsection{Example 4.} 
\begin{figure}
	\centering
	\begin{subfigure}[b]{0.45\textwidth}
		\includegraphics[width=\textwidth]{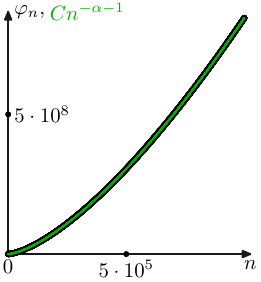}
		\caption{Exact values and first asymptotic term}
		\label{fig5a}
	\end{subfigure}
	~ 
	\begin{subfigure}[b]{0.45\textwidth}
		\includegraphics[width=\textwidth]{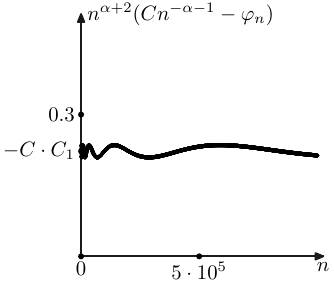}
		\caption{Difference between exact values and first asymptotic term}
		\label{fig5b}
	\end{subfigure}
	\caption{Comparison between $\vp_n$ and its asymptotics for Example 4.}\label{fig5}
\end{figure}

\begin{figure}
	\includegraphics[width=\textwidth]{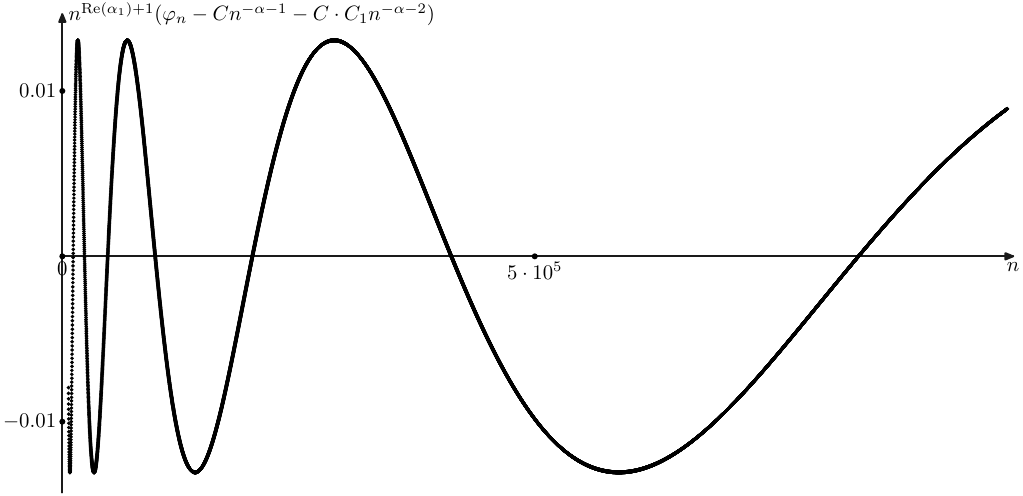}
	\caption{Oscillatory factors in the second asymptotic term for $\vp_n$, see \er{207}.}
	\label{fig6}
\end{figure}

Let us compare the above mentioned cases with a process in a random environment averaged on an incomplete interval. We take $q(z)=z$, the same as in Example 0, and fix interval $r\in(2,+\iy)$ instead of $(0,+\iy)$. Then the generating function for the ratios of rare events $\Phi(z)$ satisfies the Schr\"oder-type functional equation
$$
\int_2^{+\iy}\Phi(ze^{r(z-1)})dr=e^{-2}\Phi(z),\ \ \Phi(0)=0,\ \ \Phi'(0)=1,
$$
or
\[\lb{204}
 \int_0^{ze^{2(z-1)}}\frac{\Phi(w)}{w}dw=e^{-2}(1-z)\Phi(z),\ \ \Phi(0)=0,\ \ \Phi'(0)=1,
\]
after the change of variable similar to \er{100a}. Identity \er{204} leads to the explicit recurrence relation for determining Taylor coefficients $\vp_n$: 
$$
\vp_1=1,\ \ \ \frac{\vp_n-\vp_{n-1}}{e^2}=\frac{\vp_n}{ne^{2n}}+\frac{2(n-1)\vp_{n-1}}{1!(n-1)e^{2(n-1)}}+\frac{2^2(n-2)^2\vp_{n-2}}{2!(n-2)e^{2(n-2)}}+...+\frac{2^{n-1}\vp_{1}}{(n-1)!e^2},\ n>1.
$$ 
Again, non-recurrence formulas like \er{010} for $\vp_n$ in this example are unknown. We follow the ideas explained in \cite{K1} to find approximations of $\vp_n$ for large $n$. First of all, substituting linear combinations of $(1-z)^{\a}$ into \er{204}, we find approximation of $\Phi(z)$ at $z=1$, namely
\[\lb{205}
 \Phi(z)\simeq c\lt((1-z)^{\a}+\frac{(4\a-1)(\a+1)}{3(2\a+1)}(1-z)^{\a+1}\rt)+b(1-z)^{\a_1}+\ol{b}(1-z)^{\ol{\a_1}}+...,
\]
where 
\[\lb{206}
 \a=-2.469874943242969...,\ \ \a_1=-1.469486150158083...+{\bf i}4.38645551777719...
\]
and $\a_2=\ol{\a_1}$ are the first roots with minimal real parts of the equation
$$
 \frac{-1-\a}{e^2}=3^{\a+1}.
$$ 
Applying \er{100} and \er{101a} to \er{205}, we obtain the asymptotic expansion of $\vp_n$:
\[\lb{207}
 \vp_n\simeq C(n^{-\a-1}+C_1n^{-\a-2})+
 n^{-\Re(\a_1)-1}(B_1\cos(\Im(\a_1)\ln n)+B_2\sin(\Im(\a_1)\ln n))+...,
\]
\[\lb{207a}
 C_1:=\frac{(1+\a)(1-2\a)(2+\a)}{6(1+2\a)}.
\]
In contrast to the first Examples, we have no fast algorithms for the computation of $C$, $B_1$, and $B_2$. There is an approximate value
\[\lb{208}
 C=1.281889848285....
\]
The approximation \er{207} is quite good, see Figs. \ref{fig5} and \ref{fig6}. The oscillations in Fig. \ref{fig5b} attenuate slowly since the difference between $-\Re(\a_1)-1$ and $-\a-2$ is quite small. The difference with the standard Galton-Watson process in Example 3 is the absence of oscillations in the leading asymptotic term. However, the oscillations appear in the next term - this is the difference with the averaging over complete intervals considered in Examples 0, 1, and 2. We'll leave an intriguing question to interested readers: What happens with oscillations if we average over $(1,+\iy)$ instead of $(2,+\iy)$?

\subsection{Example 5.} As a bonus, let us consider the Galton-Watson process in random environments with the binomial offspring probability-generating function
\[\lb{300}
 P_r(z)=z(r+(1-r)z)^m,
\]
where $m\in\N$. The parameter $r$ is uniformly distributed in the interval $(0,1)$. Then the generalized Schr\"oder type integral-functional equation, see \er{004}, is
\[\lb{301}
 \int_0^1\Phi(z(r+(1-r)z)^m)dr=\frac{\Phi(z)}{m+1},\ \ \Phi(0)=0,\ \ \Phi'(0)=1.
\]
We will follow the ideas presented in Example 2 of \cite{K1}, sometimes omitting details. Changing the variable $w=z(r+(1-r)z)^m$, we rewrite \er{301} as
\[\lb{302}
 \int_{z^{m+1}}^z\Phi(w)w^{\frac1m-1}dw=\frac{m}{m+1}(1-z)z^{\frac1{m}}\Phi(z).
\]
Introducing
\[\lb{303}
 F(z)=\int_0^z\Phi(w)w^{\frac1m-1}dw,
\]
we rewrite \er{302} as
\[\lb{304}
 \frac{F(z)-F(z^{m+1})}{z-z^2}=\frac{m}{m+1}F'(z).
\]
Equations \er{303} and \er{304} leads to the recurrence formulas for determining Taylor coefficients $\vp_n$ of $\Phi(z)$:
\[\lb{304a}
 \vp_1=1,\ \ \ \vp_n=\ca 
         \frac{mn+1}{m(n-1)}\vp_{n-1},& m+1\nmid n-1,\\
         \frac{mn+1}{m(n-1)}\vp_{n-1}-\frac{(m+1)^2}{m(n-1)}\vp_{\frac{n-1}{m+1}},& otherwise
       \ac,\ n>1.
\]
Substituting {\it ansatz} $F(z)\simeq C(1-z)^{\a}$, $z\to1$ into \er{304}, we obtain the equation for $\a$:
\[\lb{304b}
 (m+1)^{\a}-1=\frac{m\a}{m+1}.
\]
The solution of \er{304b} with a minimal real part is $\a=-1$. This observation and \er{303} leads to $\Phi(z)\simeq C(1-z)^{-2}+...$ for $z\to1$ and, hence, we obtain the approximation 
\[\lb{304c}
 \vp_n\simeq Cn+....
\]
On the other hand $F(z)\sim z^{1+\frac1m}/(1+1/m)$ for $z\to0$ as it is seen from \er{303}. These two asymptotics of $F(z)$ for $z\to0$ and $z\to1$ guess us the following
definition
\[\lb{305}
 H(z)=\frac{(1-z)F(z)}{z^{1+\frac{1}{m}}}
\]
for determining $C$, which satisfies
\[\lb{305a}
 H(0)=\frac{m}{m+1},\ \ \ H(1)=C.
\]
Definition \er{305} allows us to write \er{304} as
\[\lb{306}
 -\frac{(m+1)H(z^{m+1})z^m(1-z)}{1-z^{m+1}}=m(1-z)H'(z)-H(z)
\]
or
\[\lb{307}
 -\frac{(m+1)H(z^{m+1})z^m(1-z)^{\frac{1}{m}}}{1-z^{m+1}}=m(1-z)^{\frac{1}{m}}H'(z)-(1-z)^{\frac{1}{m}-1}H(z),
\]
which is
\[\lb{308}
 -\frac{(m+1)H(z^{m+1})z^m(1-z)^{\frac{1}{m}}}{1-z^{m+1}}=(m(1-z)^{\frac{1}{m}}H(z))'.
\]
After integration \er{308} becomes
\[\lb{309}
 H(z)=\frac{(1-z)^{-\frac1m}}{m}\int_{z}^1\frac{(m+1)H(\z^{m+1})\z^m(1-\z)^{\frac{1}{m}}}{1-\z^{m+1}}d\z.
\]
Now, this is an exciting moment - using \er{309} with the integration by parts yields
\begin{multline}\lb{310}
 \int_0^w\frac{H(z)}{1-z}dz=\int_0^w\frac{(1-z)^{-1-\frac1m}dz}{m}\int_{z}^1\frac{(m+1)H(\z^{m+1})\z^m(1-\z)^{\frac{1}{m}}}{1-\z^{m+1}}d\z=\\
 (1-z)^{-\frac1m}\int_{z}^1\frac{(m+1)H(\z^{m+1})\z^m(1-\z)^{\frac{1}{m}}}{1-\z^{m+1}}d\z\lt|_{z=0}^w+\int_0^w\frac{(m+1)H(\z^{m+1})\z^m}{1-\z^{m+1}}d\z,
\end{multline}
which with \er{309} gives
\[\lb{311}
 \int_{w^{m+1}}^w\frac{H(z)}{1-z}dz=m(H(w)-H(0)).
\]
When $w\to1$, \er{311} leads to
\[\lb{312}
 H(1)\ln (m+1)=m(H(1)-H(0)),
\]
which, see \er{305a}, gives
\[\lb{313}
 C=\frac{m^2}{(m+1)(m-\ln(m+1))}.
\]
Thus, by \er{304c}, we have
\[\lb{314}
 \vp_n=\frac{m^2n}{(m+1)(m-\ln(m+1))}+...,\ \ n\to+\iy.
\]
There are still some important technical details that we should take into account to make the proof rigorous. We will skip some of them because they are partially formulated as an exercise on the site {\it math.stackexchange.com}, see \cite{K1}. We do not break the intrigue. Going back from \er{313} to \er{303}, we see that $\Phi(z)\sim C(1-z)^{-2}$ for $z\to1$. All the Taylor coefficients of $\Phi(z)$ are non-negative because they are ratios of probabilities. The non-negativity can also be obtained directly from \er{304a}. The integral $z^{-\frac1m}\int_0^zw^{\frac1m-1}\Phi(w)dw\sim C(1-z)^{-1}$ for $z\to1$. Its Taylor coefficients are also non-negative. Thus, by the Tauberian Hardy-Littlewood theorem in the Titchmarsh formulation for the series, see p. 226 in \cite{T}, we have
\[\lb{THLT}
 \sum_{j=1}^n\frac{\vp_j}{j+\frac1m}\sim Cn.
\]
Using \er{THLT} and \er{320} derived below, we obtain $\vp_n\sim Cn$.

\begin{figure}
	\centering
	\begin{subfigure}[b]{0.45\textwidth}
		\includegraphics[width=\textwidth]{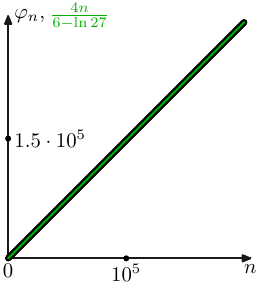}
		\caption{Exact values and first asymptotic term}
		\label{fig7a}
	\end{subfigure}
	~ 
	\begin{subfigure}[b]{0.45\textwidth}
		\includegraphics[width=\textwidth]{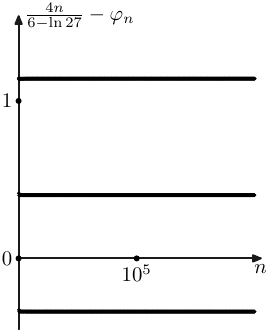}
		\caption{Difference between exact values and first asymptotic term}
		\label{fig7b}
	\end{subfigure}
	\caption{Comparison between $\vp_n$ and its asymptotics for Example 5, see \er{314}, where $r=2$.}\label{fig7}
\end{figure}

Example $m=1$ is already considered in \cite{K1}. For $m=2$, the comparison of exact values $\vp_n$, computed by \er{304a}, and its asymptotic \er{314} is provided in Fig. \ref{fig7}. The linear growth of $\vp_n$ is confirmed by its asymptotic \er{314}. The second asymptotic term is provided in Fig. \ref{fig7b}. This is a sequence of period $3$. 

\begin{figure}
	\centering
	\begin{subfigure}[b]{0.35\textwidth}
		\includegraphics[width=\textwidth]{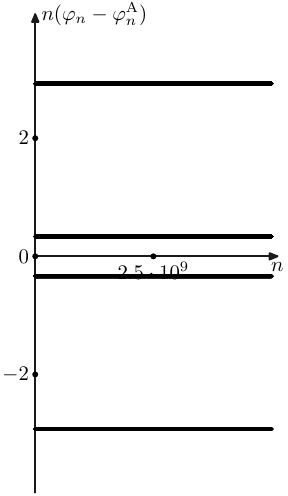}
		\caption{$m=1$}
		\label{fig8a}
	\end{subfigure}
	~ 
	\begin{subfigure}[b]{0.35\textwidth}
		\includegraphics[width=\textwidth]{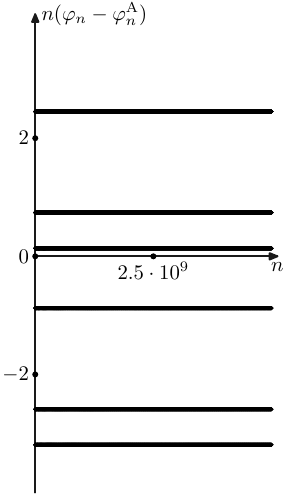}
		\caption{$m=2$}
		\label{fig8b}
	\end{subfigure}
	\hfill
	\begin{subfigure}[b]{0.35\textwidth}
		\includegraphics[width=\textwidth]{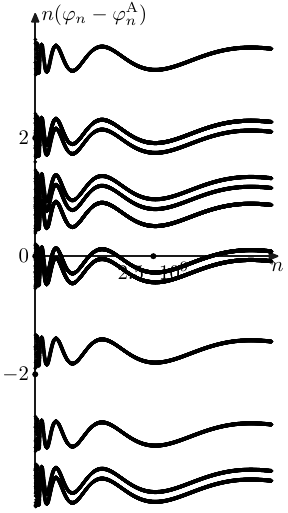}
		\caption{$m=3$}
		\label{fig8c}
	\end{subfigure}
	~ 
	\begin{subfigure}[b]{0.35\textwidth}
		\includegraphics[width=\textwidth]{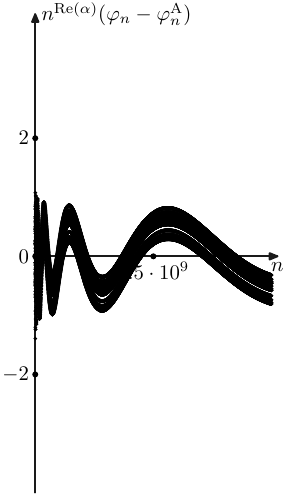}
		\caption{$m=4$, for $\a$ see \er{324}}
		\label{fig8d}
	\end{subfigure}
	\caption{The normalized reminder in \er{315} for different $m$. The  first two terms are denoted by $\vp_n^{\rm A}$.}\label{fig8}
\end{figure}

In principle, the fact that the second asymptotic term is a sequence of period $m+1$ is seen from \er{304a}. Substituting
\[\lb{315}
 \vp_n=C(n+\p_{n}+...)
\]
into \er{304a} and assuming $n\to+\iy$, we obtain
\[\lb{316}
 \p_n=\ca 
 \p_{n-1}+\frac1m,& m+1\nmid n-1,\\
 \p_{n-1}-1,& otherwise,
 \ac
\]
which means that $\p_n$ is a piece-wise linear $m+1$-periodic sequence with the slope $1/m$. To determine all the values of this sequence, it is sufficient to evaluate one of its values. Substituting Taylor expansion of $\Phi(z)$ into \er{302}, we obtain
\[\lb{318}
 \sum_{n=1}^{+\iy}\vp_n\frac{z^{(m+1)n+1}-z^n}{(n+\frac1m)(z-1)}= \frac{m}{m+1}\sum_{n=1}^{+\iy}\vp_nz^n,
\]
which is
\[\lb{319}
 \sum_{n=1}^{+\iy}\frac{\vp_nz^n}{n+\frac1m}(1+z+...+z^{mn})= \frac{m}{m+1}\sum_{n=1}^{+\iy}\vp_nz^n.
\]
Identity \er{319} yields
\[\lb{320}
 \frac{m}{m+1}\vp_{(m+1)n}=\frac{\vp_{(m+1)n}}{(m+1)n+\frac1m}+\frac{\vp_{(m+1)n-1}}{(m+1)n-1+\frac1m}+...+\frac{\vp_n}{n+\frac1m}.
\]
Substituting \er{315} into \er{320}, using $m+1$-periodicity of $\p_n$, and assuming $n\to+\iy$, we obtain
\[\lb{321}
 \frac{m}{m+1}\p_{m+1}=1-\frac1m\ln(m+1)+\frac{\p_{1}+\p_{2}+...+\p_{m+1}}{m+1}\ln(m+1),
\]
where $\sum_{j=n}^{j=(m+1)n}1/j\to\ln(m+1)$ is additionally used. Identity \er{316} along with \er{321} give
\[\lb{322}
 \p_{m+1}=\frac{m+1}{2m}\cdot\frac{(m+2)\ln(m+1)-2m}{(m+1)\ln(m+1)-m}.
\]
Results \er{T2f1} and \er{T2f2} of Theorem \ref{T2} follow from \er{315}, \er{316}, and \er{322}. Further research can be aimed at finding good estimates of the remainder in \er{315}. Here, we compute the reminder numerically, see Fig. \ref{fig8}. For $m=1,2$, see Figs. \ref{fig8a} and \ref{fig8b}, the normalized reminder tends to a periodic sequence with a discrete constant period (short-phase oscillation). For $m=3$, see Fig. \ref{fig8c}, the situation is similar to $m=1,2$. However, the convergence to the short-phase oscillation is slow, because the real part of the complex root of \er{304b} responsible for the long-phase oscillation,
\[\lb{323}
 \a=1.069829398878181...+{\bf i}5.3614900352974979...
\]
is quite close to $1$ - the value $1$ is responsible for the short-phase oscillation. Such $\a$, see \er{323}, is strictly greater than $1$, so the short-phase oscillation dominates. For $m=4$, see Fig. \ref{fig8c}, the real part of the smallest root of \er{304b} lying in the open right half-plane of $\C$
\[\lb{324}
 \a=0.870774007425338...+{\bf i}4.612162886836734...
\]
is already smaller than $1$. Thus, the long-phase oscillations become dominant. The curves slowly merge into one curve representing this long-phase oscillation, see Fig. \ref{fig8d}. The merging is slow because $\Re\a$ is close to $1$, see \er{324}. Recall also that the log-period of long-phase oscillations is determined by the imaginary part of the corresponding root, see details in \cite{K1}. For all $m\ge 4$, there are complex roots of \er{304b} with real parts less than $1$. Hence, the third term in the asymptotic will contain long-phase oscillations. This is a significant difference with the cases $m=1,2,3$. Note that there are no roots of \er{304b} with negative real parts except $\a=-1$ already considered. I am sure that the attentive reader still has important questions regarding this rich example. It will be beneficial to understand them yourself. 

\section*{Appendix} Let us start with the so-called {\it ratios of probabilities of rare events}
\[\lb{A1}
 \vp_{n,t}=\frac{\mathbb{P}(X_t=n)}{\mathbb{P}(X_t=1)}
\]
for the classical Galtot-Watson process $X_t$ with the offspring probability-generating function
\[\lb{A2}
 P(z)=p_1z+p_2z^2+p_3z^3+...,
\]
where at least one particle is survived $p_1>0$. It is well known that the characteristic function of $X_t$ has the form
\[\lb{A3}
 \underbrace{P\circ...\circ P}_{t}(z)=\mathbb{P}(X_t=1)z+\mathbb{P}(X_t=2)z^2+\mathbb{P}(X_t=3)z^3+...,
\] 
see, e.g., \cite{H}. Taking into account the obvious identity $\mathbb{P}(X_t=1)=p_1^t$, we arrive at
\[\lb{A4}
 p_1^{-t}\underbrace{P\circ...\circ P}_{t}(z)=\vp_{1,t}z+\vp_{2,t}z^2+\vp_{3,t}z^3+....
\]
Denoting
\[\lb{A5}
 \Phi(z)=\lim_{t\to+\iy}p_1^{-t}\underbrace{P\circ...\circ P}_{t}(z),\ \ \ \vp_{n}=\lim_{t\to+\iy}\vp_{n,t},
\]
we obtain
\[\lb{A6}
 \Phi(z)=\vp_1z+\vp_2z^2+\vp_3z^3+...,
\]
where, by the first identity in \er{A5}, the generating function $\Phi(z)$ satisfies the Schr\"oder functional equation with the initial conditions
\[\lb{A7}
 \Phi(P(z))=p_1\Phi(z),\ \ \ \Phi(0)=0,\ \ \ \Phi'(0)=1.
\]
For the Galton-Watson process in random environments, where $P_r(z)$ is not fixed as in \er{A2} but has a parameter $r$ which, at each time step $t$, is chosen randomly from some probabilistic space $(R,\m)$, we should write
\[\lb{A3R}
 \int_{R}...\int_{R}P_{r_1}\circ...\circ P_{r_t}(z)\m(dr_1)...\m(dr_t)=\mathbb{P}(X_t=1)z+\mathbb{P}(X_t=2)z^2+\mathbb{P}(X_t=3)z^3+...
\]
instead of \er{A3}. Going through the corresponding analogs of \er{A4}-\er{A6}, we arrive at the generalized Schr\"oder-type functional equation \er{004}, which is the analog of \er{A7} in the case of random environments. The details are provided in \cite{K1}.

For an intuitive understanding of LHS of \er{A3R}, it is enough to imagine, for example, a two-point space $(R,\m)$, then the integrals will turn into sums, and the terms will be all possible compositions of length $t$ of two functions of various orders multiplied by the probabilities of the realizations of these composition orders.

\section*{Acknowledgements} 
This paper is a contribution to the project M3 and S1 of the Collaborative Research Centre TRR 181 "Energy Transfer in Atmosphere and Ocean" funded by the Deutsche Forschungsgemeinschaft (DFG, German Research Foundation) - Projektnummer 274762653. Comments and remarks from reviewers and editors improved the quality of the manuscript. The author expresses gratitude to them.

\section*{Conflict of interest}

This study has no conflicts of interest.

\section*{Data availability statement}

Data sharing does not apply to this article as no datasets were generated or analyzed during the current study. The program code for procedures and functions considered in the article is available from the corresponding author on reasonable request.

\end{document}